\documentclass[12pt]{article}

\catcode`\@=11

\thicklines
\newskip\Einheit \Einheit=.6cm
\newcount\xcoord \newcount\ycoord
\newdimen\xdim \newdimen\ydim \newdimen\PfadD@cke \newdimen\Pfadd@cke
\def\PfadDicke#1{\PfadD@cke#1 \divide\PfadD@cke by2 
\Pfadd@cke\PfadD@cke \multiply\PfadD@cke by2}
\long\def\LOOP#1\REPEAT{\def\BODY{#1}\ITERATE}
\def\ITERATE{\BODY \let\next\ITERATE \else\let\next\relax\fi \next}
\let\REPEAT=\fi
\def\Punkt{\hbox{\raise-2pt\hbox to0pt{\hss\scriptsize$\bullet$\hss}}}

\def\DuennPunkt(#1,#2){\unskip
  \raise#2 \Einheit\hbox to0pt{\hskip#1 \Einheit
          \raise-1.5pt\hbox to0pt{\hss\tiny$\bullet$\hss}\hss}}
		  
\def\NormalPunkt(#1,#2){\unskip
  \raise#2 \Einheit\hbox to0pt{\hskip#1 \Einheit
          \raise-3pt\hbox to0pt{\hss\large$\bullet$\hss}\hss}}
\def\DickPunkt(#1,#2){\unskip
  \raise#2 \Einheit\hbox to0pt{\hskip#1 \Einheit
          \raise-4pt\hbox to0pt{\hss\Large$\bullet$\hss}\hss}}
\def\Kreis(#1,#2){\unskip
  \raise#2 \Einheit\hbox to0pt{\hskip#1 \Einheit
          \raise-4pt\hbox to0pt{\hss\Large$\circ$\hss}\hss}}
\def\Diagonale(#1,#2)#3{\unskip\leavevmode
  \xcoord#1\relax \ycoord#2\relax
      \raise\ycoord \Einheit\hbox to0pt{\hskip\xcoord \Einheit
         \unitlength\Einheit
         \line(1,1){#3}\hss}}
\def\AntiDiagonale(#1,#2)#3{\unskip\leavevmode
  \xcoord#1\relax \ycoord#2\relax \advance\xcoord by -0.05\relax
      \raise\ycoord \Einheit\hbox to0pt{\hskip\xcoord \Einheit
         \unitlength\Einheit
         \line(1,-1){#3}\hss}}
\def\Pfad(#1,#2),#3\endPfad{\unskip\leavevmode
  \xcoord#1 \ycoord#2 \thicklines\ZeichnePfad#3\endPfad\thinlines}
\def\ZeichnePfad#1{\ifx#1\endPfad\let\next\relax
  \else\let\next\ZeichnePfad
    \ifnum#1=1
      \raise\ycoord \Einheit\hbox to0pt{\hskip\xcoord \Einheit
         \vrule height\Pfadd@cke width1 \Einheit depth\Pfadd@cke\hss}%
      \advance\xcoord by 1
     \else\ifnum#1=2
      \raise\ycoord \Einheit\hbox to0pt{\hskip\xcoord \Einheit
         \unitlength\Einheit
         \line(0,1){1}\hss}
      \advance\xcoord by 0
      \advance\ycoord by 1
 \else\ifnum#1=3
      \raise\ycoord \Einheit\hbox to0pt{\hskip\xcoord \Einheit
         \unitlength\Einheit
         \line(1,1){1}\hss}
      \advance\xcoord by 1
      \advance\ycoord by 1
    \else\ifnum#1=4
      \raise\ycoord \Einheit\hbox to0pt{\hskip\xcoord \Einheit
         \unitlength\Einheit
         \line(1,-1){1}\hss}
      \advance\xcoord by 1
      \advance\ycoord by -1
   \else\ifnum#1=5
      \raise\ycoord \Einheit\hbox to0pt{\hskip\xcoord \Einheit
         \unitlength\Einheit
         \line(2,1){2}\hss}
      \advance\xcoord by 2
      \advance\ycoord by 1
	  \else\ifnum#1=6
      \raise\ycoord \Einheit\hbox to0pt{\hskip\xcoord \Einheit
         \unitlength\Einheit
         \line(2,-1){2}\hss}
      \advance\xcoord by 2
      \advance\ycoord by -1
	  \else\ifnum#1=7
      \raise\ycoord \Einheit\hbox to0pt{\hskip\xcoord \Einheit
         \unitlength\Einheit
         \line(3,1){3}\hss}
      \advance\xcoord by 3
      \advance\ycoord by 1
	  \else\ifnum#1=8
      \raise\ycoord \Einheit\hbox to0pt{\hskip\xcoord \Einheit
         \unitlength\Einheit
         \line(3,-1){3}\hss}
      \advance\xcoord by 3
      \advance\ycoord by -1
    \fi\fi\fi\fi\fi\fi\fi\fi
  \fi\next}
\def\hSSchritt{\leavevmode\raise-.4pt\hbox 
to0pt{\hss.\hss}\hskip.2\Einheit
  \raise-.4pt\hbox to0pt{\hss.\hss}\hskip.2\Einheit
  \raise-.4pt\hbox to0pt{\hss.\hss}\hskip.2\Einheit
  \raise-.4pt\hbox to0pt{\hss.\hss}\hskip.2\Einheit
  \raise-.4pt\hbox to0pt{\hss.\hss}\hskip.2\Einheit}
\def\vSSchritt{\vbox{\baselineskip.2\Einheit\lineskiplimit0pt
\hbox{.}\hbox{.}\hbox{.}\hbox{.}\hbox{.}}}
\def\DSSchritt{\leavevmode\raise-.4pt\hbox to0pt{%
  \hbox to0pt{\hss.\hss}\hskip.2\Einheit
  \raise.2\Einheit\hbox to0pt{\hss.\hss}\hskip.2\Einheit
  \raise.4\Einheit\hbox to0pt{\hss.\hss}\hskip.2\Einheit
  \raise.6\Einheit\hbox to0pt{\hss.\hss}\hskip.2\Einheit
  \raise.8\Einheit\hbox to0pt{\hss.\hss}\hss}}
\def\dSSchritt{\leavevmode\raise-.4pt\hbox to0pt{%
  \hbox to0pt{\hss.\hss}\hskip.2\Einheit
  \raise-.2\Einheit\hbox to0pt{\hss.\hss}\hskip.2\Einheit
  \raise-.4\Einheit\hbox to0pt{\hss.\hss}\hskip.2\Einheit
  \raise-.6\Einheit\hbox to0pt{\hss.\hss}\hskip.2\Einheit
  \raise-.8\Einheit\hbox to0pt{\hss.\hss}\hss}}
\def\SPfad(#1,#2),#3\endSPfad{\unskip\leavevmode
  \xcoord#1 \ycoord#2 \ZeichneSPfad#3\endSPfad}
\def\ZeichneSPfad#1{\ifx#1\endSPfad\let\next\relax
  \else\let\next\ZeichneSPfad
    \ifnum#1=1
      \raise\ycoord \Einheit\hbox to0pt{\hskip\xcoord \Einheit
         \hSSchritt\hss}%
      \advance\xcoord by 1
    \else\ifnum#1=2
      \raise\ycoord \Einheit\hbox to0pt{\hskip\xcoord \Einheit
        \hbox{\hskip-2pt \vSSchritt}\hss}%
      \advance\ycoord by 1
    \else\ifnum#1=3
      \raise\ycoord \Einheit\hbox to0pt{\hskip\xcoord \Einheit
         \DSSchritt\hss}
      \advance\xcoord by 1
      \advance\ycoord by 1
    \else\ifnum#1=4
      \raise\ycoord \Einheit\hbox to0pt{\hskip\xcoord \Einheit
         \dSSchritt\hss}
      \advance\xcoord by 1
      \advance\ycoord by -1
    \fi\fi\fi\fi
  \fi\next}
\def\Koordinatenachsen(#1,#2){\unskip
 \hbox to0pt{\hskip-.5pt\vrule height#2 \Einheit width.5pt depth1 
\Einheit}%
 \hbox to0pt{\hskip-1 \Einheit \xcoord#1 \advance\xcoord by1
    \vrule height0.25pt width\xcoord \Einheit depth0.25pt\hss}}
\def\Koordinatenachsen(#1,#2)(#3,#4){\unskip
 \hbox to0pt{\hskip-.5pt \ycoord-#4 \advance\ycoord by1
    \vrule height#2 \Einheit width.5pt depth\ycoord \Einheit}%
 \hbox to0pt{\hskip-1 \Einheit \hskip#3\Einheit 
    \xcoord#1 \advance\xcoord by1 \advance\xcoord by-#3 
    \vrule height0.25pt width\xcoord \Einheit depth0.25pt\hss}}
\def\Gitter(#1,#2){\unskip \xcoord0 \ycoord0 \leavevmode
  \LOOP\ifnum\ycoord<#2
    \loop\ifnum\xcoord<#1
      \raise\ycoord \Einheit\hbox to0pt{\hskip\xcoord 
\Einheit\Punkt\hss}%
      \advance\xcoord by1
    \repeat
    \xcoord0
    \advance\ycoord by1
  \REPEAT}
\def\Gitter(#1,#2)(#3,#4){\unskip \xcoord#3 \ycoord#4 \leavevmode
  \LOOP\ifnum\ycoord<#2
    \loop\ifnum\xcoord<#1
      \raise\ycoord \Einheit\hbox to0pt{\hskip\xcoord 
\Einheit\Punkt\hss}%
      \advance\xcoord by1
    \repeat
    \xcoord#3
    \advance\ycoord by1
  \REPEAT}
\def\Label#1#2(#3,#4){\unskip \xdim#3 \Einheit \ydim#4 \Einheit
  \def\lo{\advance\xdim by-.5 \Einheit \advance\ydim by.5 \Einheit}%
  \def\llo{\advance\xdim by-.25cm \advance\ydim by.5 \Einheit}%
  \def\loo{\advance\xdim by-.5 \Einheit \advance\ydim by.25cm}%
  \def\o{\advance\ydim by.25cm}%
  \def\ro{\advance\xdim by.5 \Einheit \advance\ydim by.5 \Einheit}%
  \def\rro{\advance\xdim by.25cm \advance\ydim by.5 \Einheit}%
  \def\roo{\advance\xdim by.5 \Einheit \advance\ydim by.25cm}%
  \def\l{\advance\xdim by-.30cm}%
  \def\r{\advance\xdim by.30cm}%
  \def\lu{\advance\xdim by-.5 \Einheit \advance\ydim by-.6 \Einheit}%
  \def\llu{\advance\xdim by-.25cm \advance\ydim by-.6 \Einheit}%
  \def\luu{\advance\xdim by-.5 \Einheit \advance\ydim by-.30cm}%
  \def\u{\advance\ydim by-.30cm}%
  \def\ru{\advance\xdim by.5 \Einheit \advance\ydim by-.6 \Einheit}%
  \def\rru{\advance\xdim by.25cm \advance\ydim by-.6 \Einheit}%
  \def\ruu{\advance\xdim by.5 \Einheit \advance\ydim by-.30cm}%
  #1\raise\ydim\hbox to0pt{\hskip\xdim
     \vbox to0pt{\vss\hbox to0pt{\hss$#2$\hss}\vss}\hss}%
}
\catcode`\@=12

\usepackage{amsmath,amsthm}
\usepackage{color}
\usepackage{xspace}
\usepackage[colorlinks=true,
linkcolor=green,
filecolor=brown,
citecolor=green]{hyperref}

\def\v{\vert}
\def\lr{LRmax\xspace}
\def\a{\ensuremath{\mathcal A}\xspace}

\def\c{\ensuremath{\mathcal C}\xspace}

\def\si{\ensuremath{\sigma}}
\def\ok{$3\underline{5}241$OK\xspace}

\begin{document}
\newtheorem{lemma}{Lemma}
\newtheorem{theorem}{Theorem}
\newtheorem{prop}{Proposition}
\newtheorem{cor}{Corollary}
\vspace*{-10mm}
\begin{center}
{\Large                        
A Combinatorial Interpretation of the Eigensequence for Composition
}

\vspace{10mm}
DAVID CALLAN  \\
Department of Statistics  \\
\vspace*{-2mm}
University of Wisconsin-Madison  \\
\vspace*{-2mm}
Medical Science Center \\
\vspace*{-2mm}
1300 University Ave  \\
\vspace*{-2mm}
Madison, WI \ 53706-1532  \\
{\bf callan@stat.wisc.edu}  \\
\vspace{5mm}

July 20, 2005
\end{center}

\vspace{3mm}
\begin{center}
   \textbf{Abstract}
\end{center}
    
 The monic sequence that shifts left under convolution with itself is the 
Catalan numbers with 130+ combinatorial interpretations. Here we 
establish a combinatorial interpretation for the monic sequence that shifts 
left under \emph{composition}: it counts permutations that contain a 3241 pattern only 
as part of a 35241 pattern. We give two recurrences, 
the first allowing relatively fast computation,
the second similar to one for the Catalan numbers. 
Among the $4\times 4!=96$ similarly restricted patterns involving 4 
letters (such as $4\underline{2}31$: a 431 pattern only occurs as part 
of a 4231), four different counting sequences arise: 64 give the Catalan 
numbers, 16 give the Bell numbers, 12 give sequence 
\htmladdnormallink{A051295}{http://www.research.att.com:80/cgi-bin/access.cgi/as/njas/sequences/eisA.cgi?Anum=A051295}
in OEIS, and 4 give a new sequence with an explicit formula.   

\vspace{10mm}

{\Large \textbf{1 \quad Introduction}  }

The composition of two sequences $(a_{n})_{n\ge1}$ and $(b_{n})_{n\ge1}$ 
is $(c_{n})_{n\ge1}$ defined by $C(x)=A(B(x))$ where $A,B,C$ are the 
respective generating functions, $A(x)=\sum_{n\ge 1}a_{n}x^{n}$ and so 
on.
A sequence is monic if its first term is 1. There is a 
unique monic sequence $(b_{n})_{n\ge1}$  whose composition with itself 
is equal to its left shift, $(b_{2},b_{3},\ldots)$. This sequence 
$(b_{n})_{n\ge1}$ begins $1,1,2,6,23,104,531,\ldots$ and is called an 
eigensequence for composition \cite{canonical} .

Consider a permutation on a set of positive integers as a list (or 
word of distinct letters). A subpermutation is a subword (letters in 
same order but not not necessarily contiguous). Thus 253 is a 
subpermutation of 21534. A factor is a subpermutation in which the 
letters are contiguous. A standard permutation is one on an initial 
segment of the positive integers. The reduced form, 
reduce($\pi$), of a permutation $\pi$ on an arbitrary set of positive 
integers is the standard permutation obtained by replacing its 
smallest entry by 1, next smallest by 2, and so on. Thus 
reduce$(352)=231$. The length $\v \pi \v$ of  a 
permutation $\pi$ is simply the number of letters in it. Given a permutation $\pi$ 
and a standard permutation $\rho$ of weakly smaller length---a 
``pattern''---a subpermutation of $\pi$ whose reduced form is $\rho$ is 
said to be an instance of the pattern $\rho$ in $\pi$. Thus 352 forms 
a 231 pattern in 35124. A 321-avoiding permutation, for example, is 
one containing no instance of a 321 pattern. The number of 321-avoiding 
permutations is the Catalan number $C_{n}$ \cite{ec2}.

Let $\a_{n}$ denote the set of permutations on $[n]$ in which the 
pattern 3241 only occurs as part of a 35241 pattern 
($3\underline{5}241$OK permutations for short). This curious pattern 
restriction arises in connection with 2-stack-sortable permutations, 
which can be characterized as permutations that are both \ok and 
2341-avoiding \cite{dulucq1}. With $a_{n}=\vert \a_{n} \vert$, the 
first several terms of $(a_{n})_{n\ge 0}$ coincide with those of 
$(b_{n})_{n\ge1}$  above. This suggests that $a_{n}=b_{n+1}$ and the 
main objective of this paper is to prove bijectively that this is so (\S 3 
and \S 4). In \S 2, we consider the structure of a \ok permutation 
and deduce two recurrences for $a_{n}$. A bonus section (\S 5) 
classifies all similarly restricted patterns involving 4 letters by 
their counting sequences.

\vspace*{6mm}

{\Large \textbf{2 \quad Structure and Recurrences for 
$\mathbf{3\underline{5}241}$OK Permutations}  }

Every permutation $\pi$ on $[n]$ has the form 
$\pi=m_{1}L_{1}m_{2}L_{2}\ldots m_{r}L_{r}$ where 
$m_{1}<m_{2}<\ldots<m_{r}=n$ are the left-to-right maxima (\lr for 
short) of $\pi$. With this notation, the following characterization of \ok permutations is 
easy to verify.
\begin{theorem}
    A permutation $\pi$ on $[n]$ is \emph{\ok} if and only if
    
    $($\mbox{\,i}\,$)$ $L_{1}<L_{2}<\ldots<L_{r}$ in the sense that $u\in 
    L_{i},\,v\in L_{j}$ with $i<j $ implies $u<v$, and 
    
    $($ii\,$)$ each $L_{i}$ is \ok. \qed
\end{theorem}
Recall $\a_{n}$ is the set of all \ok permutations on $[n]$.
Let $\a_{n,k}$ 
denote the subset of $\a_{n}$ with first entry $k$, $\c_{n}$ the subset 
of $\a_{n}$ with first two entries increasing ($n\ge 2,\ 
\c_{1}:=\a_{1}=\{(1)\}$). Set $a_{n}=\v\a_{n}\v,\ a_{n,k}=\v \a_{n,k} 
\v,\ c_{n}=\v \c_{n} \v$.
\begin{theorem}
    $a_{n}$ is given by the recurrence relations $a_{0}=c_{1}=1$ and
    
     \mbox{\ $($i}\,$)\qquad a_{n}=\sum_{i=0}^{n-1}a_{i}c_{n-i}$ 
     \hspace*{34mm} $n \ge 1$
     
     \mbox{\,$($ii}\,$)\qquad c_{n}=\sum_{i=0}^{n-1} i a_{n-1,i}$  
      \hspace*{33mm} $n \ge 2$
      
       $($\mbox{iii}\,$)\qquad a_{n,k}=
       \begin{cases}
	   \sum_{i=0}^{k-1}a_{i}\sum_{j=k-i}^{n-1-i}a_{n-1-i,j} &
	   \quad 1 \le k \le n-1 \\
	   a_{n-1} & \quad k=n
	   \end{cases}$
\end{theorem}

\noindent \textbf{Proof} \quad (i) counts $\a_{n}$ by $i=\v L_{1} \v$ 
since $m_{1}L_{1}K \longrightarrow 
$ \big(reduce($L_{1}$),\:reduce($m_{1}K$)\big) 
is a bijection to $\a_{i} \times \c_{n-i}$. (ii) counts $\c_{n}$ by 
second entry, say $i+1$, because $L_{1}=\emptyset$ and, with 
$m_{2}=i+1,\ m_{1}m_{2}L_{2}\ldots m_{r}L_{r} \longrightarrow \ 
\big(m_{1},\:$reduce($m_{2}L_{2}\ldots m_{r}L_{r}$)\big) is a bijection to  
$[i] \times \a_{n-1,i}$. (iii) counts $\a_{n,k}\ (1\le k \le n-1)$ 
by $i=\v L_{1} \v$ and $j=$ first entry of 
reduce($m_{2}L_{2}\ldots$) since $j \ge k-i,\, m_{1}=k$ and 
$m_{1}L_{1}m_{2}L_{2}\ldots \longrightarrow \big( 
L_{1},\:$reduce$(m_{2}L_{2}\ldots)\big)$ is a bijection to $\a_{i} 
\times \a_{n-1-i,j}$.   \qed

There is a more elegant (but less computationally efficient) recurrence 
involving a sum over compositions. Recall that a composition 
$\mathbf{c}$ of $n$ is a list $\mathbf{c}=(c_{1},c_{2},\ldots, c_{r})$ 
of positive integers whose sum is $n$ and if 
$\mathbf{c}=(c_{i})_{i=1}^{r}$ and $\mathbf{d}=(d_{i})_{i=1}^{r}$ 
are same-length compositions of $n$, then $\mathbf{d}$ dominates 
$\mathbf{c}\ (\mathbf{d} \succeq \mathbf{c})$ if $d_{1}+\ldots+d_{i} 
\ge c_{1}+\ldots+c_{i}$ for $i=1,2,\ldots,r-1$ (of course, equality 
holds for $i=r$). Let $\c_{n}$ denote the set of all compositions of 
$n$ (there are $2^{n-1}$ of them).
A permutation $\pi=m_{1}L_{1}m_{2}L_{2}\ldots m_{r}L_{r}$ (recall the 
$m_{i}$ are the \lr entries) determines two same-length compositions of $[n]$: 
$\mathbf{c}=(c_{i})_{i=1}^{r}$ with $c_{i}=\v m_{i}L_{i} \v =1+\v 
L_{i} \v$ and $\mathbf{d}=(d_{i})_{i=1}^{r}$  with 
$d_{i}=m_{i}-m_{i-1}\ (m_{0}:=0)$ and necessarily $\mathbf{d} \succeq 
\mathbf{c}$ (or else a left-to-right max would occur among the 
$L_{i}$). Summing over $\mathbf{c}$, Theorem 1 yields
\begin{theorem}
    \[
     \hspace*{20mm} a_{n}=\sum_{\mathbf{c}\in \c_{n}} \v \{\mathbf{d}\in \c_{n}\,:\, 
    \mathbf{d} \succeq \mathbf{c}\}\v\  
    a_{c_{1}-1}a_{c_{2}-1}\cdots a_{c_{r}-1} \hspace*{20mm} n\ge 1 
    \] \qed
\end{theorem}
Omitting the first factor in the summand, the recurrence 
$a_{0}=1,\ a_{n} = \sum_{\mathbf{c}\in \c_{n}} 
a_{c_{1}-1}a_{c_{2}-1}\cdots a_{c_{r}-1}$ ($n\ge 1$) is well known to 
generate the Catalan numbers (for example, count Dyck paths by the 
locations of their returns to ground level).

\vspace*{6mm}

{\Large \textbf{3 \quad Preparing for the Main Bijection}  }

The generating function $B(x)=\sum_{n\ge 1}b_{n}x^{n}$ for the 
``shifts left under composition'' sequence of \S 1 is characterized 
by $B(B(x))=\frac{B(x)}{x}-1$. 
If $(b_{n})$ is the counting sequence for some species, say 
B-structures, then $[x^{n}]B(B(x))$ is the number of pairs 
$(X,\mathbf{Y})$ where $X$ is a B-structure of unspecified size $k,\ 
1 \le k \le n$ and $\mathbf{Y} $ is a $k$-list of B-structures of 
total size $n$ (a $k$-list is simply a list with $k$ entries). On the other hand,
$[x^{n}](\frac{B(x)}{x}-1)=b_{n+1}$. Hence, to show that 
$b_{n+1}=a_{n}$, our main objective, we need a bijection from the 
set $\a_{n}$ of \ok permutations on $[n]$ to pairs 
$(\rho,\mathbf{v})$ where $\rho \in \a_{k-1}$ for some $1\le k \le n$ 
and $\mathbf{v}$ is a $k$-list of \ok permutations (possibly empty) of 
total length $n-k$.

Indeed, given $\pi \in \a_{n}$, the position of $n$ in $\pi$ 
determines $k$ and the right factor of $\pi$ starting at $n$ 
determines $\rho$ as follows. It is convenient to define the LIT 
(longest increasing terminal) entries of a permutation $\pi$ on $[n]$ 
to be $k,k+1,\ldots,n$ where $k$ is the smallest integer such that 
$k,k+1,\ldots,n$ appear in that order in $\pi$. Note that the LIT 
entries form a terminal segment of the \lr entries. For example, the 
LIT entries are double underlined and the remaining \lr entries are 
single underlined in 
$\underline{2}\,1\,\underline{4}\,\underline{\underline{7}}\,6\,5 \,
\underline{\underline{8}}\,\underline{\underline{9}}\,3$. The LIT 
entries of a permutation on arbitrary positive integers are defined 
analogously. Now decompose 
$\pi$ as $\si n \tau$.
\begin{theorem}
    A permutation $\si n \tau$ on $[n]$ is \emph{\ok} if and only if 
    $\si,\tau$ are both \emph{\ok} and each entry of $\si$ that exceeds the 
    smallest entry of $\tau$ is an LIT entry of $\si$.
\end{theorem}

\noindent \textbf{Proof} \quad $(\Rightarrow)$ \quad Suppose $b \in 
\si$ exceeds $c \in \tau$ yet $b$ is not an LIT entry of $\si$. Then 
either (i) there exists $a \in \si$ preceding $b$ with $a>b$, and $abnc$ is an 
offending 3241 pattern, or (ii) there exist $u,v$ following $b$ in 
$\si$ with 
$u>v>b$ and this time $uvnc$ is an offending 3241.

\hspace*{10mm}  $(\Leftarrow)$ \quad routine. \qed

Set $k=\v n\tau\v$, the length of the right 
factor of $\pi$ starting at $n$, and
set $\rho=$ reduce($\tau) \in \a_{k-1}$, capturing the underlying 
permutation of $\tau$. The support of $\tau$ can be captured by 
placing, for each $c \in$\,support($\tau)$, an asterisk (or star) right before 
the smallest entry of $\si$ that exceeds $c$ (necessarily an LIT 
entry in view of Theorem 4) or after $\max(\si)$ 
if there is no such entry. No information is lost if $\si$ is 
then reduced (keeping the stars in place) to produce $\si^{*}$ 
say, because only the LIT entries of $\si$ will be affected and the 
stars determine both the original LIT entries of $\si$ and the 
support of $\tau$. For example, $
\pi=(2,8,3,1,11,4,6,5,13,7,15,9,10,14,12) \in \a_{15}$ gives $k=5,\ 
\tau=(9,10,14,12)$ and yields the pair $\rho=(1,2,4,3),\ 
\si^{*}=(2, 8, 3, 1,*,*, 9, 4, 6, 5, *, 10, *, 7)$, and $\pi$ can be 
recovered from $\rho$ and $\si^{*}$.
The length of $\si^{*}$ (disregarding the stars) coincides with 
the total length $n-k$ of the desired list of \ok permutations. This 
reduces the problem to giving a bijection from \ok permutations on 
$[n]$ with 
$k-1$ stars distributed arbitrarily, just before LIT entries or 
just after the maximum entry $n$, to $k$-lists of (possibly empty) \ok 
permutations of total length $n$.

Let $X_{n,k}$ denote 
the subset of the preceding ``starred'' permutations on $[n]$ with (i) 
no stars immediately following or preceding the max $n$, and (ii) at most one 
star preceding each non-max LIT entry. In other words, $X_{n,k}$ is the set 
of \ok permutations on $[n]$ with $k-1$ distinct LIT entries other 
than $n$ preceded (or marked) by a star. We now show that the following result 
suffices to construct the desired bijection.
\begin{theorem}
    With $X_{n,k}$ as just defined, there is a bijection from $X_{n,k}$ 
    to the set of $k$-lists of \emph{nonempty} \ok permutations of 
    total length $n$.
\end{theorem}
Consider a permutation with (possibly) 
multiple stars in the allowed locations, for example, with $k=8$ and 
hence
7 stars and with $a,b$ denoting non-max LIT entries,
\[
\ldots***a\ldots*b\ldots*n**\ldots
\]
We will show how to represent this as an element of $X_{n,j}$ for 
some $j$ together with a bit sequence of $j$ 0s and $k-j$ 1s. The 
preceding Theorem then converts the element of $X_{n,j}$ to a $j$-list 
of nonempty \ok permutations of total length $n$, while the 1s in the 
bit sequence specify the locations for empty permutations and we have the 
desired $k$-list of permutations of total length $n$.

First, collapse all contiguous stars to a single star and delete the 
stars (if any) surrounding $n$. The example yields
\[
\ldots***a\ldots*b\ldots*n**\ldots \quad \rightarrow \quad
\ldots*a\ldots*b\ldots n\ldots\quad \in X_{n,j}
\]
with $j=3$ and $j-1=2$ stars.
Second, replace each star immediately preceding a non-max LIT entry 
and $n$ itself by 1, all other stars by 0, and then suppress the 
permutation entries to get a bit sequence $j$ 1s and $k-j$ 0s; the 
example yields
\[
\ldots***a\ldots*b\ldots*n**\ldots \quad \rightarrow \quad
(0\,0\,\overset{*a}{1}\,\overset{*b}{1}\,0\,\overset{n}{1}\,0\,0).
\]
Clearly, the original ``multiple-starred'' permutation can be uniquely 
recovered from the element of $X_{n,j}$ and the bit sequence, and we 
are all done as soon as we establish Theorem 5. In fact, we will prove 
Theorem 5 with a specific type of bijection that permits a further 
reduction in the problem.

Define the \lr factors of a permutation $m_{1}L_{1}m_{2}L_{2}\ldots 
m_{r}L_{r}$ to be $m_{1}L_{1},m_{2}L_{2},$ $\ldots, m_{r}L_{r}$. Then there 
is a refined form of Theorem 5 that goes as follows.
\begin{theorem}
    Let $X_{n,k}$ denote the set of \emph{\ok} permutations on $[n]$
    with $k-1$ of the non-max LIT entries marked. Let $Y_{n,k}$ 
    denote the set of $k$-lists of nonempty \ok permutations of total 
    length $n$. There is a bijection from $X_{n,k}$ to $Y_{n,k}$, 
    $\pi \rightarrow \mathbf{v}$, of the following type: it rearranges 
    the \lr factors of $\pi$, then splits them into a $k$-list of 
    permutations and, finally, reduces each one to yield $\mathbf{v}$.
\end{theorem}
It suffices to present this bijection for the special case where
$L_{1},L_{2},\ldots,L_{r}$ are all increasing lists, equivalently, 
since $L_{1}<L_{2}<\ldots <L_{r}$ for a \ok permutation, the entire 
list $L_{1}L_{2}\ldots L_{r}$ is increasing. (In the general case, sort 
each $L_{i}$ and then apply the special case bijection except, before 
reducing, replace each sorted $L_{i}$ by the original list $L_{i}$.)

Note that permutations $\pi=m_{1}L_{1}m_{2}L_{2}\ldots 
m_{r}L_{r}$ for which the full list $L_{1}L_{2}\ldots L_{r}$ is increasing are 
precisely the 321-avoiding permutations (and automatically \ok). In summary, we have shown in 
this section that the entire problem boils down to
\begin{theorem}
    There is a bijection of the following specific type from $321$-avoiding permutations $\pi$ on $[n]$ with 
    $k-1$ of the non-max LIT entries marked to $k$-lists $\mathbf{v}$ of nonempty 
    $321$-avoiding permutations of total 
    length $n$: 
    it rearranges 
    the \lr factors of $\pi$, then splits them into a $k$-list of 
    permutations and, finally, reduces each one to yield $\mathbf{v}$.
\end{theorem}

\vspace*{6mm}

{\Large \textbf{4 \quad The Main Bijection}  }

We now prove Theorem 7 and use 
\[ 
\underline{3}\,1\,\underline{5}\,2\,\underline{8}\,
4\,6\,\underline{12}\,7\,\underline{15}\,9\,\underline{17}\,10\,11\,
\underline{20}\,\underline{\underline{25}}\,
\overset{*}{\underline{\underline{26}}}\,13\,\underline{\underline{27}}\, 
\overset{*}{\underline{\underline{28}}}\,14\,
\overset{*}{\underline{\underline{29}}}\,16\,
\underline{\underline{30}}\,18\,19\,21\,22\,23\,24
\]
with $n=30$ and $k=4$ as a working example (as before, the LIT entries are 
double underlined and the other \lr entries are single underlined). 
First, we split the permutation into panes each consisting of one or 
more \lr factors. To do so, use a moving window initially covering a 
right factor of the permutation and consisting of $k$ panes starting, 
respectively, at the very first LIT entry and the first LIT entry 
following each starred LIT entry. (Since the window gradually moves 
left, dropping the last pane and sometimes acquiring a new first 
pane, it is drawn below the permutation for clarity.)
\vspace*{2mm}
\Einheit=0.5cm
\[
\Label\o{ \textrm{{3}}}(-14.5,1)
\Label\o{ \textrm{{1}}}(-13.5,1)
\Label\o{ \textrm{{5}}}(-12.5,1)
\Label\o{ \textrm{{2}}}(-11.5,1)
\Label\o{ \textrm{{8}}}(-10.5,1)
\Label\o{ \textrm{{4}}}(-9.5,1)
\Label\o{ \textrm{{6}}}(-8.5,1)
\Label\o{ \textrm{{12}}}(-7.5,1)
\Label\o{ \textrm{{7}}}(-6.5,1)
\Label\o{ \textrm{{15}}}(-5.5,1)
\Label\o{ \textrm{{9}}}(-4.5,1)
\Label\o{ \textrm{{17}}}(-3.5,1)
\Label\o{ \textrm{{10}}}(-2.5,1)
\Label\o{ \textrm{{11}}}(-1.5,1)
\Label\o{ \textrm{{20}}}(-.5,1)
\Label\o{ \textrm{{25}}}(.5,1)
\Label\o{ \textrm{{26}}}(1.5,1)
\Label\o{ \textrm{{13}}}(2.5,1)
\Label\o{ \textrm{{27}}}(3.5,1)
\Label\o{ \textrm{{28}}}(4.5,1)
\Label\o{ \textrm{{14}}}(5.5,1)
\Label\o{ \textrm{{29}}}(6.5,1)
\Label\o{ \textrm{{16}}}(7.5,1)
\Label\o{ \textrm{{30}}}(8.5,1)
\Label\o{ \textrm{{18}}}(9.5,1)
\Label\o{ \textrm{{19}}}(10.5,1)
\Label\o{ \textrm{{21}}}(11.5,1)
\Label\o{ \textrm{{22}}}(12.5,1)
\Label\o{ \textrm{{23}}}(13.5,1)
\Label\o{ \textrm{{24}}}(14.5,1)
\Pfad(0,1),111111111111111\endPfad
\Pfad(0,0),111111111111111\endPfad
\Pfad(0,0),2\endPfad
\Pfad(3,0),2\endPfad
\Pfad(6,0),2\endPfad
\Pfad(8,0),2\endPfad
\Pfad(15,0),2\endPfad
\SPfad(0,1),2\endSPfad
\SPfad(3,1),2\endSPfad
\SPfad(6,1),2\endSPfad
\SPfad(8,1),2\endSPfad
\SPfad(15,1),2\endSPfad
\]
\vspace*{-15mm}
\begin{center}
The Initial Window
\end{center}
Proceed as follows until the entire permutation is covered by 
nonoverlapping panes. (It's just possible that the initial window 
covers the entire permutation but this happens only for the identity 
permutation. In this case associate the symbol $\emptyset$ with the 
window and proceed to the second step below.) Find the maximum non-LRmax entry in all but the 
last pane of the current window (the max of the empty set is 
considered to be $-\infty$) and denote this number $m$. Look for 
\lr entries not yet empaned that exceed $m$. There are 2 cases.

\noindent Case 1.\quad There is no such \lr entry. In this case associate 
$\emptyset$ with the current window, and delete the last pane to obtain 
the next window.

\noindent Case 2.\quad There exist \lr entries not yet empaned that 
exceed $m$. In this case associate the smallest such entry, $M$ say, with the 
current window and shift the window left by prepending a pane 
starting at $M$ and deleting the last pane.

The resulting panes and the successive windows, each with an 
associated \lr or $\emptyset$ symbol, are shown.

\Einheit=0.5cm
\[
\Label\o{ \textrm{{3}}}(-14.5,9)
\Label\o{ \textrm{{1}}}(-13.5,9)
\Label\o{ \textrm{{5}}}(-12.5,9)
\Label\o{ \textrm{{2}}}(-11.5,9)
\Label\o{ \textrm{{8}}}(-10.5,9)
\Label\o{ \textrm{{4}}}(-9.5,9)
\Label\o{ \textrm{{6}}}(-8.5,9)
\Label\o{ \textrm{{12}}}(-7.5,9)
\Label\o{ \textrm{{7}}}(-6.5,9)
\Label\o{ \textrm{{15}}}(-5.5,9)
\Label\o{ \textrm{{9}}}(-4.5,9)
\Label\o{ \textrm{{17}}}(-3.5,9)
\Label\o{ \textrm{{10}}}(-2.5,9)
\Label\o{ \textrm{{11}}}(-1.5,9)
\Label\o{ \textrm{{20}}}(-.5,9)
\Label\o{ \textrm{{25}}}(.5,9)
\Label\o{ \textrm{{26}}}(1.5,9)
\Label\o{ \textrm{{13}}}(2.5,9)
\Label\o{ \textrm{{27}}}(3.5,9)
\Label\o{ \textrm{{28}}}(4.5,9)
\Label\o{ \textrm{{14}}}(5.5,9)
\Label\o{ \textrm{{29}}}(6.5,9)
\Label\o{ \textrm{{16}}}(7.5,9)
\Label\o{ \textrm{{30}}}(8.5,9)
\Label\o{ \textrm{{18}}}(9.5,9)
\Label\o{ \textrm{{19}}}(10.5,9)
\Label\o{ \textrm{{21}}}(11.5,9)
\Label\o{ \textrm{{22}}}(12.5,9)
\Label\o{ \textrm{{23}}}(13.5,9)
\Label\o{ \textrm{{24}}}(14.5,9)
\Pfad(0,9),111111111111111\endPfad
\Pfad(0,8),111111111111111\endPfad
\Pfad(0,8),2\endPfad
\Pfad(3,8),2\endPfad
\Pfad(6,8),2\endPfad
\Pfad(8,8),2\endPfad
\Pfad(15,8),2\endPfad
\Label\r{ \rightarrow 17}(15.5,8.5)
\Pfad(-4,8),111111111111\endPfad
\Pfad(-4,7),111111111111\endPfad
\Pfad(-4,7),2\endPfad
\Pfad(0,7),2\endPfad
\Pfad(3,7),2\endPfad
\Pfad(6,7),2\endPfad
\Pfad(8,7),2\endPfad
\Label\r{ \rightarrow 15}(8.5,7.5)
\Pfad(-6,7),111111111111\endPfad
\Pfad(-6,6),111111111111\endPfad
\Pfad(-6,6),2\endPfad
\Pfad(-4,6),2\endPfad
\Pfad(0,6),2\endPfad
\Pfad(3,6),2\endPfad
\Pfad(6,6),2\endPfad
\Label\r{ \rightarrow \emptyset}(6.5,6.5)
\Pfad(-6,6),111111111\endPfad
\Pfad(-6,5),111111111\endPfad
\Pfad(-6,5),2\endPfad
\Pfad(-4,5),2\endPfad
\Pfad(0,5),2\endPfad
\Pfad(3,5),2\endPfad
\Label\r{ \rightarrow 12}(3.5,5.5)
\Pfad(-8,5),11111111\endPfad
\Pfad(-8,4),11111111\endPfad
\Pfad(-8,4),2\endPfad
\Pfad(-6,4),2\endPfad
\Pfad(-4,4),2\endPfad
\Pfad(0,4),2\endPfad
\Label\r{ \rightarrow \emptyset}(0.5,4.5)
\Pfad(-8,4),1111\endPfad
\Pfad(-8,3),1111\endPfad
\Pfad(-8,3),2\endPfad
\Pfad(-6,3),2\endPfad
\Pfad(-4,3),2\endPfad
\Label\r{ \rightarrow 8}(-3.5,3.5)
\Pfad(-11,3),11111\endPfad
\Pfad(-11,2),11111\endPfad
\Pfad(-11,2),2\endPfad
\Pfad(-8,2),2\endPfad
\Pfad(-6,2),2\endPfad
\Label\r{ \rightarrow \emptyset}(-5.5,2.5)
\Pfad(-11,2),111\endPfad
\Pfad(-11,1),111\endPfad
\Pfad(-11,1),2\endPfad
\Pfad(-8,1),2\endPfad
\Label\r{ \rightarrow 3}(-7.5,1.5)
\Pfad(-15,1),1111\endPfad
\Pfad(-15,0),1111\endPfad
\Pfad(-15,0),2\endPfad
\Pfad(-11,0),2\endPfad
\Label\r{ \rightarrow \emptyset}(-10.5,0.5)
\SPfad(-15,1),222222222\endSPfad
\SPfad(-11,3),2222222\endSPfad
\SPfad(-8,5),22222\endSPfad
\SPfad(-6,7),222\endSPfad
\SPfad(-4,8),22\endSPfad
\SPfad(0,9),2\endSPfad
\SPfad(3,9),2\endSPfad
\SPfad(6,9),2\endSPfad
\SPfad(8,9),2\endSPfad
\SPfad(15,9),2\endSPfad
\]
\begin{center}
The Moving Window
\end{center}
Second, rearrange the panes using a $k$-row array formed from a list 
consisting of (i) the 
windows' associated \lr entries and $\emptyset$s, and (ii) the first 
entries of the $k$ panes of the initial window:
\[
\emptyset\ 3\ \emptyset\ 8\ \emptyset\ 12\ \emptyset\ 15\ 17\ 25\ 27\ 29\ 30
\]
Insert the list entries right to left into array entries lower right 
to upper left moving up each column in turn EXCEPT, when a $\emptyset$ 
is inserted, that row is ``not accepting any more entries'' and is 
henceforth skipped over. The resulting array is
\[
\begin{array}{ccccc}
          &           & \emptyset &     12     & 25  \\
          &           &           & \emptyset  & 27  \\
\emptyset &  3        &     8     &     15     & 29  \\
          &           & \emptyset &     17     & 30
\end{array}
\]
Now form a list of $k$ permutations by concatenating the panes 
initiated by the \lr entries in each row:
\[
12\:7\:25\:26\:13,\quad 27\:28\:14,\quad 
3\:1\:5\:2\:8\:4\:6\:15\:9\:29\:16,\quad 
17\:10\:11\:20\:30\:18\:19\:21\:22\:23\:24
\]
Finally, reduce each one to get the desired $k$-list of 321-avoiding 
permutations:
\[
2\:1\:4\:5\:3,\quad 2\:3\:1,\quad 
3\:1\:5\:2\:7\:4\:6\:9\:8\:11\:10,\quad 
3\:1\:2\:6\:11\:4\:5\:7\:8\:9\:10.
\]
Note that LIT entries in these reduced permutations correspond to LIT 
entries in the original as do \lr entries.

To reverse the process, suppose given a list 
$(\pi_{1},\pi_{2},\ldots,\pi_{k})$ of $k$ 321-avoiding permutations of 
total length $n$ and let us use the previous example with $k=4$ and 
$n=30$. Double underline the LIT entries in each $\pi_{i}$ and single 
underline the remaining \lr entries. Form a pane from the first 
LIT entry to the end of each permutation.

\Einheit=0.5cm
\[
\Label\o{ \underline{2}}(-16.5,0)
\Label\o{1}(-15.5,0)
\Label\o{ \underline{\underline{4}}}(-14.5,0)
\Label\o{ \pi_{1}}(-14.5,1.5)
\Label\o{ \underline{\underline{5}}}(-13.5,0)
\Label\o{3}(-12.5,0)
\Pfad(-15,0),111\endPfad
\Pfad(-15,1),111\endPfad
\Pfad(-15,0),2\endPfad
\Pfad(-12,0),2\endPfad
\Label\o{,}(-11.5,-0.5)
\Label\o{ \underline{\underline{2}}}(-10.5,0)
\Label\o{ \underline{\underline{3}}}(-9.5,0)
\Label\o{ \pi_{2}}(-9.5,1.5)
\Label\o{1}(-8.5,0)
\Pfad(-11,0),111\endPfad
\Pfad(-11,1),111\endPfad
\Pfad(-11,0),2\endPfad
\Pfad(-8,0),2\endPfad
\Label\o{,}(-7.5,-0.5)
\Label\o{ \underline{3}}(-6.5,0)
\Label\o{1}(-5.5,0)
\Label\o{ \underline{5}}(-4.5,0)
\Label\o{2}(-3.5,0)
\Label\o{ \underline{7}}(-2.5,0)
\Label\o{4}(-1.5,0)
\Label\o{\pi_{3}}(-1.5,1.5)
\Label\o{6}(-0.5,0)
\Label\o{ \underline{9}}(0.5,0)
\Label\o{8}(1.5,0)
\Label\o{ \underline{\underline{11}}}(2.5,0)
\Label\o{10}(3.5,0)
\Pfad(2,0),11\endPfad
\Pfad(2,1),11\endPfad
\Pfad(2,0),2\endPfad
\Pfad(4,0),2\endPfad
\Label\o{,}(4.5,-0.5)
\Label\o{\underline{3}}(5.5,0)
\Label\o{1}(6.5,0)
\Label\o{2}(7.5,0)
\Label\o{\underline{6}}(8.5,0)
\Label\o{ \underline{\underline{11}}}(9.5,0)
\Label\o{4}(10.5,0)
\Label\o{\pi_{4}}(10.5,1.5)
\Label\o{5}(11.5,0)
\Label\o{7}(12.5,0)
\Label\o{8}(13.5,0)
\Label\o{9}(14.5,0)
\Label\o{10}(15.5,0)
\Pfad(9,0),1111111\endPfad
\Pfad(9,1),1111111\endPfad
\Pfad(9,0),2\endPfad
\Pfad(16,0),2\endPfad
\]
\begin{center}
The Initial Panes
\end{center}
Place $n,n-1,n-2,\ldots$ below the LIT entries working right to left.
\newpage 
\Einheit=0.5cm
\[
\Label\o{ \underline{2}}(-16.5,0)
\Label\o{1}(-15.5,0)
\Label\o{ \underline{\underline{4}}}(-14.5,0)
\Label\o{ \underline{\underline{5}}}(-13.5,0)
\Label\o{ 25}(-14.5,-1)
\Label\o{ 26}(-13.5,-1)
\Label\o{3}(-12.5,0)
\Pfad(-15,-1),111\endPfad
\Pfad(-15,0),111\endPfad
\Pfad(-15,1),111\endPfad
\Pfad(-15,-1),22\endPfad
\Pfad(-12,-1),22\endPfad
\Label\o{,}(-11.5,-0.5)
\Label\o{ \underline{\underline{2}}}(-10.5,0)
\Label\o{ \underline{\underline{3}}}(-9.5,0)
\Label\o{ 27}(-10.5,-1)
\Label\o{28}(-9.5,-1)
\Label\o{1}(-8.5,0)
\Pfad(-11,-1),111\endPfad
\Pfad(-11,0),111\endPfad
\Pfad(-11,1),111\endPfad
\Pfad(-11,-1),22\endPfad
\Pfad(-8,-1),22\endPfad
\Label\o{,}(-7.5,-0.5)
\Label\o{ \underline{3}}(-6.5,0)
\Label\o{1}(-5.5,0)
\Label\o{ \underline{5}}(-4.5,0)
\Label\o{2}(-3.5,0)
\Label\o{ \underline{7}}(-2.5,0)
\Label\o{4}(-1.5,0)
\Label\o{6}(-0.5,0)
\Label\o{ \underline{9}}(0.5,0)
\Label\o{8}(1.5,0)
\Label\o{ \underline{\underline{11}}}(2.5,0)
\Label\o{29}(2.5,-1)
\Label\o{10}(3.5,0)
\Pfad(2,-1),11\endPfad
\Pfad(2,0),11\endPfad
\Pfad(2,1),11\endPfad
\Pfad(2,-1),22\endPfad
\Pfad(4,-1),22\endPfad
\Label\o{,}(4.5,-0.5)
\Label\o{\underline{3}}(5.5,0)
\Label\o{1}(6.5,0)
\Label\o{2}(7.5,0)
\Label\o{\underline{6}}(8.5,0)
\Label\o{ \underline{\underline{11}}}(9.5,0)
\Label\o{30}(9.5,-1)
\Label\o{4}(10.5,0)
\Label\o{5}(11.5,0)
\Label\o{7}(12.5,0)
\Label\o{8}(13.5,0)
\Label\o{9}(14.5,0)
\Label\o{10}(15.5,0)
\Pfad(9,-1),1111111\endPfad
\Pfad(9,0),1111111\endPfad
\Pfad(9,1),1111111\endPfad
\Pfad(9,-1),22\endPfad
\Pfad(16,-1),22\endPfad
\]
\begin{center}
Place LIT Entries
\end{center}
Set $b=$ largest integer in $[n]$ not yet placed; here $b=24$. Start 
with the last permutation in the list. Take the largest entry that is 
blank below (here, 10), place $b$ in the blank, decrement $b$ and 
repeat as long as these largest entries are (i) empaned or (ii) underlined. 
Then add a pane that begins at the leftmost newly added entry unless 
it is already empaned.

\Einheit=0.5cm
\[
\Label\o{ \underline{2}}(-16.5,0)
\Label\o{1}(-15.5,0)
\Label\o{ \underline{\underline{4}}}(-14.5,0)
\Label\o{ \underline{\underline{5}}}(-13.5,0)
\Label\o{ 25}(-14.5,-1)
\Label\o{ 26}(-13.5,-1)
\Label\o{3}(-12.5,0)
\Pfad(-15,-1),111\endPfad
\Pfad(-15,0),111\endPfad
\Pfad(-15,1),111\endPfad
\Pfad(-15,-1),22\endPfad
\Pfad(-12,-1),22\endPfad
\Label\o{,}(-11.5,-0.5)
\Label\o{ \underline{\underline{2}}}(-10.5,0)
\Label\o{ \underline{\underline{3}}}(-9.5,0)
\Label\o{ 27}(-10.5,-1)
\Label\o{28}(-9.5,-1)
\Label\o{1}(-8.5,0)
\Pfad(-11,-1),111\endPfad
\Pfad(-11,0),111\endPfad
\Pfad(-11,1),111\endPfad
\Pfad(-11,-1),22\endPfad
\Pfad(-8,-1),22\endPfad
\Label\o{,}(-7.5,-0.5)
\Label\o{ \underline{3}}(-6.5,0)
\Label\o{1}(-5.5,0)
\Label\o{ \underline{5}}(-4.5,0)
\Label\o{2}(-3.5,0)
\Label\o{ \underline{7}}(-2.5,0)
\Label\o{4}(-1.5,0)
\Label\o{6}(-0.5,0)
\Label\o{ \underline{9}}(0.5,0)
\Label\o{8}(1.5,0)
\Label\o{ \underline{\underline{11}}}(2.5,0)
\Label\o{29}(2.5,-1)
\Label\o{10}(3.5,0)
\Pfad(2,-1),11\endPfad
\Pfad(2,0),11\endPfad
\Pfad(2,1),11\endPfad
\Pfad(2,-1),22\endPfad
\Pfad(4,-1),22\endPfad
\Label\o{,}(4.5,-0.5)
\Label\o{\underline{3}}(5.5,0)
\Label\o{17}(5.5,-1)
\Label\o{1}(6.5,0)
\Label\o{2}(7.5,0)
\Label\o{\underline{6}}(8.5,0)
\Label\o{20}(8.5,-1)
\Label\o{ \underline{\underline{11}}}(9.5,0)
\Label\o{30}(9.5,-1)
\Label\o{4}(10.5,0)
\Label\o{5}(11.5,0)
\Label\o{7}(12.5,0)
\Label\o{8}(13.5,0)
\Label\o{9}(14.5,0)
\Label\o{10}(15.5,0)
\Label\o{18}(10.5,-1)
\Label\o{19}(11.5,-1)
\Label\o{21}(12.5,-1)
\Label\o{22}(13.5,-1)
\Label\o{23}(14.5,-1)
\Label\o{24}(15.5,-1)
\Pfad(5,-1),11111111111\endPfad
\Pfad(5,0),11111111111\endPfad
\Pfad(5,1),11111111111\endPfad
\Pfad(5,-1),22\endPfad
\Pfad(9,-1),22\endPfad
\Pfad(16,-1),22\endPfad
\]
\begin{center}
Work Last Permutation
\end{center}
Decrement $b$ (now $b=16$) and proceed to the next permutation to the left 
in the list (considering the last to be the left neighbor of the first). Repeat until $b=0$ and no blanks 
remain. The second step gives

\Einheit=0.5cm
\[
\Label\o{ \underline{2}}(-16.5,0)
\Label\o{1}(-15.5,0)
\Label\o{ \underline{\underline{4}}}(-14.5,0)
\Label\o{ \underline{\underline{5}}}(-13.5,0)
\Label\o{ 25}(-14.5,-1)
\Label\o{ 26}(-13.5,-1)
\Label\o{3}(-12.5,0)
\Pfad(-15,-1),111\endPfad
\Pfad(-15,0),111\endPfad
\Pfad(-15,1),111\endPfad
\Pfad(-15,-1),22\endPfad
\Pfad(-12,-1),22\endPfad
\Label\o{,}(-11.5,-0.5)
\Label\o{ \underline{\underline{2}}}(-10.5,0)
\Label\o{ \underline{\underline{3}}}(-9.5,0)
\Label\o{ 27}(-10.5,-1)
\Label\o{28}(-9.5,-1)
\Label\o{1}(-8.5,0)
\Pfad(-11,-1),111\endPfad
\Pfad(-11,0),111\endPfad
\Pfad(-11,1),111\endPfad
\Pfad(-11,-1),22\endPfad
\Pfad(-8,-1),22\endPfad
\Label\o{,}(-7.5,-0.5)
\Label\o{ \underline{3}}(-6.5,0)
\Label\o{1}(-5.5,0)
\Label\o{ \underline{5}}(-4.5,0)
\Label\o{2}(-3.5,0)
\Label\o{ \underline{7}}(-2.5,0)
\Label\o{4}(-1.5,0)
\Label\o{6}(-0.5,0)
\Label\o{ \underline{9}}(0.5,0)
\Label\o{8}(1.5,0)
\Label\o{ \underline{\underline{11}}}(2.5,0)
\Label\o{15}(0.5,-1)
\Label\o{29}(2.5,-1)
\Label\o{16}(3.5,-1)
\Label\o{10}(3.5,0)
\Pfad(0,-1),1111\endPfad
\Pfad(0,0),1111\endPfad
\Pfad(0,1),1111\endPfad
\Pfad(0,-1),22\endPfad
\Pfad(2,-1),22\endPfad
\Pfad(4,-1),22\endPfad
\Label\o{,}(4.5,-0.5)
\Label\o{\underline{3}}(5.5,0)
\Label\o{17}(5.5,-1)
\Label\o{1}(6.5,0)
\Label\o{2}(7.5,0)
\Label\o{\underline{6}}(8.5,0)
\Label\o{20}(8.5,-1)
\Label\o{ \underline{\underline{11}}}(9.5,0)
\Label\o{30}(9.5,-1)
\Label\o{4}(10.5,0)
\Label\o{5}(11.5,0)
\Label\o{7}(12.5,0)
\Label\o{8}(13.5,0)
\Label\o{9}(14.5,0)
\Label\o{10}(15.5,0)
\Label\o{18}(10.5,-1)
\Label\o{19}(11.5,-1)
\Label\o{21}(12.5,-1)
\Label\o{22}(13.5,-1)
\Label\o{23}(14.5,-1)
\Label\o{24}(15.5,-1)
\Pfad(5,-1),11111111111\endPfad
\Pfad(5,0),11111111111\endPfad
\Pfad(5,1),11111111111\endPfad
\Pfad(5,-1),22\endPfad
\Pfad(9,-1),22\endPfad
\Pfad(16,-1),22\endPfad
\]
\begin{center}
Work Penultimate Permutation
\end{center}
and the final result (omitting the no-longer-needed underlines) is

\Einheit=0.5cm
\[
\Label\o{2}(-16.5,0)
\Label\o{1}(-15.5,0)
\Label\o{12}(-16.5,-1)
\Label\o{7}(-15.5,-1)
\Label\o{ 4}(-14.5,0)
\Label\o{5}(-13.5,0)
\Label\o{ 25}(-14.5,-1)
\Label\o{ 26}(-13.5,-1)
\Label\o{3}(-12.5,0)
\Label\o{13}(-12.5,-1)
\Pfad(-17,-1),11111\endPfad
\Pfad(-17,0),11111\endPfad
\Pfad(-17,1),11111\endPfad
\Pfad(-17,-1),22\endPfad
\Pfad(-15,-1),22\endPfad
\Pfad(-12,-1),22\endPfad
\Label\o{,}(-11.5,-0.5)
\Label\o{2}(-10.5,0)
\Label\o{3}(-9.5,0)
\Label\o{ 27}(-10.5,-1)
\Label\o{28}(-9.5,-1)
\Label\o{1}(-8.5,0)
\Label\o{14}(-8.5,-1)
\Pfad(-11,-1),111\endPfad
\Pfad(-11,0),111\endPfad
\Pfad(-11,1),111\endPfad
\Pfad(-11,-1),22\endPfad
\Pfad(-8,-1),22\endPfad
\Label\o{,}(-7.5,-0.5)
\Label\o{3}(-6.5,0)
\Label\o{1}(-5.5,0)
\Label\o{5}(-4.5,0)
\Label\o{2}(-3.5,0)
\Label\o{7}(-2.5,0)
\Label\o{4}(-1.5,0)
\Label\o{6}(-0.5,0)
\Label\o{3}(-6.5,-1)
\Label\o{1}(-5.5,-1)
\Label\o{5}(-4.5,-1)
\Label\o{2}(-3.5,-1)
\Label\o{8}(-2.5,-1)
\Label\o{4}(-1.5,-1)
\Label\o{6}(-0.5,-1)
\Label\o{9}(0.5,0)
\Label\o{8}(1.5,0)
\Label\o{9}(1.5,-1)
\Label\o{11}(2.5,0)
\Label\o{15}(0.5,-1)
\Label\o{29}(2.5,-1)
\Label\o{16}(3.5,-1)
\Label\o{10}(3.5,0)
\Pfad(-7,-1),11111111111\endPfad
\Pfad(-7,0),11111111111\endPfad
\Pfad(-7,1),11111111111\endPfad
\Pfad(-7,-1),22\endPfad
\Pfad(-3,-1),22\endPfad
\Pfad(0,-1),22\endPfad
\Pfad(2,-1),22\endPfad
\Pfad(4,-1),22\endPfad
\Label\o{,}(4.5,-0.5)
\Label\o{3}(5.5,0)
\Label\o{17}(5.5,-1)
\Label\o{1}(6.5,0)
\Label\o{2}(7.5,0)
\Label\o{10}(6.5,-1)
\Label\o{11}(7.5,-1)
\Label\o{6}(8.5,0)
\Label\o{20}(8.5,-1)
\Label\o{11}(9.5,0)
\Label\o{30}(9.5,-1)
\Label\o{4}(10.5,0)
\Label\o{5}(11.5,0)
\Label\o{7}(12.5,0)
\Label\o{8}(13.5,0)
\Label\o{9}(14.5,0)
\Label\o{10}(15.5,0)
\Label\o{18}(10.5,-1)
\Label\o{19}(11.5,-1)
\Label\o{21}(12.5,-1)
\Label\o{22}(13.5,-1)
\Label\o{23}(14.5,-1)
\Label\o{24}(15.5,-1)
\Pfad(5,-1),11111111111\endPfad
\Pfad(5,0),11111111111\endPfad
\Pfad(5,1),11111111111\endPfad
\Pfad(5,-1),22\endPfad
\Pfad(9,-1),22\endPfad
\Pfad(16,-1),22\endPfad
\]
\begin{center}
Final Result
\end{center}
The marked LIT entries are now retrieved as the max entry in the last 
pane of all but the last permutation, here 26,\,28,\,29. Arrange 
the panes in increasing order of first entries to retrieve the original 
starred permutation.

\vspace*{6mm}

{\Large \textbf{5 \quad Underlined 4-Patterns}  }

There are $4\times 4!=96$ similarly restricted patterns involving 4 
letters (let's call them underlined 4-patterns), for instance, $4\underline{2}31$: a 431 pattern only occurs as part 
of a 4231. We say that a permutation meeting this pattern restriction satisfies $4\underline{2}31$ or
is $4\underline{2}31$OK. Four different counting sequences arise. We will show that 64 give the Catalan 
numbers, 16 give the Bell numbers, 12 give sequence 
\htmladdnormallink{A051295}{http://www.research.att.com:80/cgi-bin/access.cgi/as/njas/sequences/eisA.cgi?Anum=A051295}
in OEIS, and 4 give the sequence with $n$th term 
$=(n-1)!+\sum_{k=0}^{n-2}\sum_{\substack{i,j\ge 0 \\
i+j \le k}}k^{\overline{i}}\,(n-2-k)^{\underline{j}}$ using 
rising/falling factorial notation.

The complement of a permutation $\pi$ on $[n]$ is $n+1-\pi$ (entrywise). The complement, reverse and inverse of an underlined
4-pattern are defined in the obvious way so that a permutation satisfies a given underlined 4-pattern if and only if its complement 
(resp. reverse, resp. inverse) satisfy the  complement 
(resp. reverse, resp. inverse) of the underlined 4-pattern. Under the action 
of the 8-element abelian group generated by the
complement, reverse and inverse operations, the 96 underlined 4-patterns 
are partitioned into equivalence classes, the members of each one all 
having the same counting sequence.

In fact, 64 underlined 4-patterns entail avoidance of the associated 
3-pattern altogether. For each of the 6 patterns of length 3, the 
permutations avoiding it are counted by the Catalan numbers \cite{ec2}. 
Hence these 64 are  counted by Catalan numbers
\htmladdnormallink{ 
A000108}{http://www.research.att.com:80/cgi-bin/access.cgi/as/njas/sequences/eisA.cgi?Anum= A000108}. (Incidentally, they 
split into 5 equivalence classes of size 8 and 6 of size 4.) The 
remaining 32---the nontrivial underlined 4-patterns---split into 5 
classes as shown in the Table (the top row contains convenient class 
representatives).

\[
\begin{array}{||c||c||c||c||c||}\hline
    3\,2\,\underline{4}\,1 &   3\,1\,\underline{4}\,2   & \underline{1}\,3\,4\,2   & \underline{1}\,3\,2\,4   & 3\,2\,1\,\underline{4}  \\ \hline
    1\,3\,4\,\underline{2} &  \underline{3}\,1\,4\,2    & \underline{1}\,4\,2\,3    & 1\,3\,2\,\underline{4}  & \underline{4}\,1\,2\,3  \\
    1\,\underline{4}\,2\,3 &    3\,1\,4\,\underline{2}  &  2\,3\,1\,\underline{4}   & 4\,2\,3\,\underline{1}  & \underline{1}\,4\,3\,2  \\
     2\,3\,\underline{1}\,4 & 3\,\underline{1}\,4\,2    & 2\,4\,3\,\underline{1}    & \underline{4}\,2\,3\,1  & 2\,3\,4\,\underline{1}  \\
     \underline{2}\,4\,3\,1 &   2\,4\,\underline{1}\,3  & 3\,1\,2\,\underline{4}    &   &    \\
     \underline{3}\,1\,2\,4 &  \underline{2}\,4\,1\,3   & 3\,2\,4\,\underline{1}    &   &   \\
     4\,\underline{1}\,3\,2 &   2\,4\,1\,\underline{3}  & \underline{4}\,1\,3\,2    &   &   \\
       4\,2\,1\,\underline{3} &  2\,\underline{4}\,1\,3 & \underline{4}\,2\,1\,3    &   &  \\\hline\hline
       \textrm{Bell} & \textrm{Bell} & \textrm{\textrm{A051295}} &  A051295 & \textrm{new} \\ \hline
\end{array}
\]
\begin{center}
The 5 equivalence classes of nontrivial underlined 4-patterns \\ and their counting sequences
\end{center}
It remains to verify the counting sequences.
\vspace*{2mm}

\noindent {\large \textbf{Bell Sequence}}\quad Consider two canonical ways to represent a 
set partition of $[n]$.
\vspace*{-2mm}
\begin{description}
    \item[Canonical Increasing]  largest entry first in each block, 
    rest of block increasing,  blocks 
    arranged in increasing order of first element, as in 412-6-735.
    \item[Canonical Decreasing]  each block decreasing, blocks 
    arranged in increasing order of first element, as in 421-6-753.
\end{description}
\begin{theorem}
    Splitting a permutation on $[n]$ into its \lr factors 
    $m_{1}L_{1},\ldots,m_{r}L_{r}$  is a 
    bijection from 
    
    $($i\,$)\   32\underline{4}1$\emph{OK} permutations on $[n]$ to set 
    partitions of $[n]$ in canonical increasing form.
    
   $ ($ii\,$)\ 31\underline{4}2$\emph{OK} permutations on $[n]$ to set 
    partitions of $[n]$ in canonical decreasing form.
\end{theorem}
\textbf{Proof}\quad (i) If $\pi$ is $32\underline{4}1$OK yet 
some $L_{i}$ fails to be increasing, then there exist $u$ preceding $v$ 
in $L_{i}$ with $u>v$ and $m_{i}uv$ is an offending 321 pattern. 
Conversely, if each $L_{i}$ is increasing, then a 321 must have its 
``2'' and ``1'' in different blocks and $\pi$ is $\pi$ is $32\underline{4}1$OK.

(ii) If $\pi$ is $31\underline{4}2$OK yet 
some \lr factor $m_{i}L_{i}$ fails to be decreasing, then there exist $u$ preceding $v$ 
in $L_{i}$ with $u<v$ and $m_{i}uv$ is an offending 312 pattern. 
Conversely, if the \lr factors of $\pi$ are all decreasing, then a 312 must have its 
``1'' and ``2'' in different blocks and $\pi$ is
$31\underline{4}2$OK. \qed

Since set partitions are counted by the Bell numbers
\htmladdnormallink{A000110}{http://www.research.att.com:80/cgi-bin/access.cgi/as/njas/sequences/eisA.cgi?Anum=A000110},
this verifies the first 2 columns in the Table.
\vspace*{2mm}

\noindent {\large \textbf{A051295 Sequence}}\quad Let us count 
$\underline{1}324$OK permutations on $[n]$. We claim each such 
permutation $\pi$ has the form
$
\pi_{1}1\pi_{2}
$ where $\pi_{1}$ is a $\underline{1}324$OK permutation on 
$[n-k+2,n]$ with $k$ the position of 1 in $\pi$, and $\pi_{2}$ is an 
arbitrary permutation on $[2,n-k+1]$. To see the claim, observe that 
each such $\pi$ is $\underline{1}324$OK. Conversely, if $\pi$ 
is $\underline{1}324$OK, all entries 
preceding 1 in $\pi$ necessarily exceed all entries following 1 or 
else there is a subpermutation $u1v$ with $u<v$ forming an offending 
324 pattern. While $\pi_{1}$ must be $\underline{1}324$OK, the 
entry 1 relieves $\pi_{2}$ of any obligation, and the claim follows.
This decomposition implies that the number $u_{n}$ of $\underline{1}324$OK permutations on $[n]$ satisfies 
$u_{n}=\sum_{k=1}^{n}u_{k-1}(n-k)!$, and hence $(u_{n})$ is sequence 
\htmladdnormallink{A051295}{http://www.research.att.com:80/cgi-bin/access.cgi/as/njas/sequences/eisA.cgi?Anum=A051295}
in OEIS.

Next, we count $\underline{1}342$OK permutations on $[n]$ by 
position $k$ of 1. If 1 occurs in first position, the rest of the permutation 
is arbitrary. Otherwise, we claim the permutation has the form (with 1 
in position $k\ge 2$)
\[
a_{1}\ldots a_{k-1}\,1\,\pi[2,a_{k-1}-1]\,\pi[a_{k-1}+1,a_{k-2}-1]\,\ldots 
\,\pi[a_{2}+1,a_{1}-1]\,
\pi[a_{1}+1,n]
\] 
where $a_{1}>a_{2}>\ldots >a_{k-1}$ and $\pi[b,c]$ denotes an 
arbitrary permutation on the interval of integers $[b,c]$ (understood 
to be the empty permutation if $b>c$).

Such permutations are clearly $\underline{1}342$OK. Conversely, 
given a $\underline{1}342$OK permutation, the entries 
$a_{1},\ldots,a_{k-1}$ preceding 
1 are necessarily decreasing left to right for otherwise there is a 
subpermutation $a\,b\,1$ with $a<b$. This is a $342$ pattern with no hope 
of the required ``1''. And for $2\le i \le k-1$, all elements of
$[a_{i}+1,a_{i-1}-1]$ precede all elements of 
$[a_{i-1}+1,a_{i-2}-1] \sqcup \ldots \sqcup [a_{1}+1,n] $ else some 
$b>a_{i-1}$ precedes some $c<a_{i-1}$ in the permutation and 
$a_{i-1}\,b\,c$ is an offending 342 pattern because all entries preceding 
$a_{i-1}$ exceed  $a_{i-1}$. The claim follows.

With $u_{n,k}$ the number of $\underline{1}342$OK permutations 
on $[n]$ for which 1 occurs in position $k$,
this decomposition implies the generating function identity for  
each $k\ge 1$
\[
\sum_{n\ge 0}u_{n,k}x^{n}=x^{k}\big( \sum_{n\ge 0}n!x^{n}\big)^{k}
\]
and so $(u_{n,k})$ forms sequence
\htmladdnormallink{A084938}{http://www.research.att.com:80/cgi-bin/access.cgi/as/njas/sequences/eisA.cgi?Anum=A084938}
in OEIS with row sums
\htmladdnormallink{A051295}{http://www.research.att.com:80/cgi-bin/access.cgi/as/njas/sequences/eisA.cgi?Anum=A051295}.
This verifies columns 3 and 4 of the Table.

The fact that $\underline{1}324$ and $\underline{1}342$ are Wilf-equivalent
(that is, have the same counting sequence) despite being in different 
equivalence classes can be explained by a simple bijection from 
$\underline{1}324$OK permutations to $\underline{1}342$OK 
permutations. Factor a $\underline{1}324$OK permutation $\pi$ as
\[
m_{1}\,L_{1}\,m_{2}\,L_{2}\,\ldots \,m_{r}\,L_{r}
\]
where $m_{1},m_{2},\ldots $ are now the left-to-right \emph{minima} of 
$\pi$. Then reassemble as
\[
m_{1}\,m_{2}\,\ldots\, m_{r}\,L_{r}\,L_{r-1}\,\ldots \, L_{1}.
\]
This works because $L_{i}$ is a permutation on $[m_{i}+1,m_{i-1}-1]$ 
for $1\le i \le r$ ($m_{0}:=n+1$) for otherwise there is an 
entry $b$ following $m_{i}$ in $\pi$ with $b>m_{i-1}$ and 
$m_{i-1}\,m_{i}\,b$ is an offending 324 pattern.

\noindent {\large \textbf{New Sequence}}\quad The 
$321\underline{4}$OK permutations on $[n]$ can be counted 
directly by position of $n$. If $n$ occurs in last position, the rest 
of the permutation is unrestricted---$(n-1)!$ choices. If $n$ occurs 
in position $n-1$ and the last entry is $i$, then $i+1,i+2,\ldots,n$ must 
occur in that order else $i$ terminates an offending 321. Thus the permutation 
is determined by the positions of 
$1,2,\ldots,i-1$---$(n-2)^{\underline{i-1}}$ choices where 
$j^{\underline{i}}=j(j-1)\ldots(j-i+1)$ is the falling factorial.

Now suppose $n$ occurs in position $k\le n-2$. The subpermutation 
following $n$ must be increasing and hence has the form $1 
\le a<y_{1}<y_{2}<\ldots<y_{n-k-2}<b<n$. Fix $a$ and $b$ and let 
$(x_{i})_{i=1}^{r}=(a,b)\backslash (y_{i})_{i=1}^{n-k-2}$ entailing
$r=(b-a-1)-(n-k-2)=b-t$ with $t=a+n-k-1$. The subpermutation preceding $n$ has no 
restriction on the entries $<a$ and placing these entries gives 
$(k-1)^{\underline{a-1}}$ choices. However, its entries $>b$ must be 
increasing left to right (else there exist $u>v>b$ with $u$ 
preceding $v$ and $u\,v\,b$ is an offending 321) and must all lie to 
the right of the $x_{i}$s (else some $u>b$ precedes some $x_{i}$ and 
$u>b>x_{i}>a$ implies $u\,x_{i}\,a$ is an offending 321). Thus the 
subpermutation preceding $n$ and involving entries  $>a$ yields only 
$r!$ choices: arrange the $x_{i}$s. All told, the desired count is 
(with $t:=a+n-k-1$ and $n\ge 1$)
\[
(n-1)!+\sum_{i=1}^{n-1}(n-2)^{\underline{i-1}} +\sum_{k=1}^{n-2} 
\sum_{a=1}^{k} 
\sum_{b=t}^{n-1}(k-1)^{\underline{a-1}}\binom{b-a-1}{n-k-2} 
(b-t)!  = 
\]
\[
(n-1)!+\sum_{k=0}^{n-2}\sum_{\substack{i,j\ge 0 \\
i+j \le k}}k^{\overline{i}}\,(n-2-k)^{\underline{j}}.
\]
This sequence begins $(1, 1, 2, 5, 15, 55, 248, 1357, 8809, 
66323, 568238,\ldots)_{n \ge 0}$ and is entrywise $\ge$ the next largest counting 
sequence
\htmladdnormallink{A051295}{http://www.research.att.com:80/cgi-bin/access.cgi/as/njas/sequences/eisA.cgi?Anum=A051295}
 $=(1,1,2,5,15,54,235,1237\ldots)$; both seem to be $\sim 
 (n-1)!$ asymptotically.\qed

A pattern restriction similar to one of our underlined 4-patterns 
arises in a recent paper \cite{patience} on ``Patience Sorting'' a deck of cards into 
piles, namely, the restriction that a 342 pattern in which the ``4'' and ``2'' 
are contiguous occurs only as part of a 3142 pattern. It is 
almost immediate, however, that the permutations satisfying this seemingly 
weaker restriction coincide with our  $3\underline{1}42$OK 
permutations (and hence are counted by the Bell numbers).

\end{document}